# A simple model on streamflow management with a dynamic risk measure

Hidekazu Yoshioka[1,*], Yumi Yoshioka[1]

[1] Graduate School of Natural Science and Technology, Shimane University, Nishikawatsu-cho 1060, Matsue, 690-8504, Japan

* Corresponding author: yoshih@life.shimane-u.ac.jp

**Abstract**

We present an exactly-solvable risk-minimizing stochastic differential game for flood management in rivers. The streamflow dynamics follow stochastic differential equations driven by a Lévy process. An entropic dynamic risk measure is employed to evaluate a flood risk under model uncertainty. The problem is solved via a Hamilton–Jacobi–Bellman–Isaacs equation. We explicitly derive an optimal flood mitigation policy along with its existence criteria and the worst-case probability measure. A backward stochastic differential representation as an alternative formulation is also presented.

**Keywords**

River management under uncertainty; Levy process; Zero-sum differential game; Hamilton–Jacobi–Bellman–Isaacs equation

# 1. Introduction

River is a dynamic environment driven by stochastic streamflow. Stochastic differential equations (SDEs) of the Ornstein–Uhlenbeck types driven by Gaussian [1] and Poisson noises [2-3] have been efficient models representing streamflow time series. Recently, a tempered stable Ornstein–Uhlenbeck model [4] has been proposed as a more realistic alternative to represent both continuous and jump noises in the streamflow time series [5].

Modeling streamflow itself is not a goal but rather a pathway toward streamflow management as dynamic optimization problems like dam and hydropower management [6-7]. Stochastic control based on the dynamic programming with SDEs [8] is a versatile tool for dynamic optimization under uncertainty. Especially, a risk-minimizing problem where the decision-maker is not fully confident about his/her reference model is a game-type control problem with a dynamic risk measure [9]. In river environmental management, it is often difficult to accurately identify the target model due to limitations such as technical constraints and lack of knowledge [10-11]. The dynamic risk measure is a mathematically rigorous tool harmonizing with this situation; however, such an approach to streamflow management does not exist to the best of our knowledge.

Solving a stochastic differential game reduces to finding an appropriate solution to a corresponding Hamilton–Jacobi–Bellman–Isaacs (HJBI) equation [8]. This kind of nonlinear partial integro-differential equation can be solved explicitly only in limited cases [12-15]. The use of such explicit solutions is of course limited; however, they provide clear and useful insights into more complicated problems.

We present an exactly-solvable streamflow management problem subject to stochastic streamflow driven by a pure jump Lévy process [16]. This is a streamflow diversion problem to minimize both flooding and diversion effort (**Figure 1**). The performance functional to be minimized by optimizing the diversion effort is set as a dynamic risk measure of an entropic type so that the deviation between the true and believed model is penalized [8]. The resulting HJBI equation is non-linear, non-local, and high-dimensional. Nevertheless, we show that it is possible to derive its explicit solution along with its existence criteria.

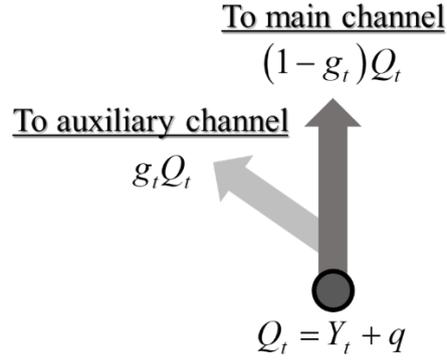

**Figure 1.** Conceptual diagram of the model.

## 2. Problem formulation

### 2.1 System dynamics

We consider a streamflow management under uncertainty as a simplified flood mitigation problem. There are the two channels, main and auxiliary channels (**Figure 1**). The goal here is to simultaneously find the optimal diversion policy and the risk measure under model uncertainty.

We assume that the streamflow $Q = (Q_t)_{t \geq 0}$ at a point in a main channel follows a Levy-driven Ornstein–Uhlenbeck model:

$$\mathrm{d}Y_t = -\lambda Y_t \mathrm{d}t + \int_0^\infty z N(\mathrm{d}z, \mathrm{d}t), \quad t > 0, \quad Y_0 \geq 0 \tag{1}$$

with $Q_t = Y_t + q$, where $q \geq 0$ is the baseflow and $Y = (Y_t)_{t \geq 0}$ is the fluctuation,

$\lambda > 0$ is the recession rate, $N$ is the Poisson measure of a Levy process having positive jumps (subordinator) and the Levy measure $v$. We assume $\int_0^\infty \min\{1, z\} v(\mathrm{d}z) < +\infty$ as in the usual setting [16]. Set the compensated Poisson measure $\tilde{N}(\mathrm{d}z, \mathrm{d}t) = N(\mathrm{d}z, \mathrm{d}t) - v(\mathrm{d}z)$. A natural filtration generated by $N$ is $\mathbb{F} = (\mathbb{F}_t)_{t \geq 0}$.

At each time $t$, the ratio $g_t \in \omega \subseteq [0,1]$ of $Q_t$ ($\omega$: a compact set) flowing into the main channel is diverted into the auxiliary channel with the instantaneous effort $n^{-1} g_t^n (Y_t + q)$, $n > 1$. This means a gradual increase of the effort with respect to $g_t$. The cumulative water volume flowing into the main channel $X = (X_t)_{t \geq 0}$ and cumulative diversion effort $U = (U_t)_{t \geq 0}$ satisfy

$$\mathrm{d}X_t = (1 - g_t)(Y_t + q) \mathrm{d}t, \quad t > 0, \quad X_0 = 0 \tag{2}$$

and

$$\mathrm{d}U_t = n^{-1} g_t^n (Y_t + q) \mathrm{d}t, \quad t > 0, \quad U_0 = 0. \tag{3}$$

The diversion effort has a specific functional shape to formulate a tractable model. Problems with general parameterizations would require some numerical methods.

Other than the streamflow dynamics, a martingale transformation between equivalent probability measures representing model uncertainty, is formulated [8]. We consider a problem where the measure $v$ is the source of uncertainty. Set the reference and distorted measures as $\mathbb{P}$ and $\mathbb{Q} \in P$, respectively, where $P$ is the collection of probability measures absolutely continuous with respect to $\mathbb{P}$. Set the Radon–Nikodým derivative $\frac{\mathrm{d}\mathbb{Q}}{\mathrm{d}\mathbb{P}}$. Consider the measure transform $M = (M_t)_{t \geq 0}$ with $\frac{\mathrm{d}\mathbb{Q}}{\mathrm{d}\mathbb{P}} = M_T$ for a terminal time $T > 0$. Its governing SDE is

$$dM_t = M_{t-}\int_0^\infty \theta(t,z)\tilde{N}(dz,dt), \quad t>0, \quad M_0 = 1 \tag{4}$$

with $\theta:[0,T]\times\mathbb{R}_+ \to (-1,+\infty)$. An element $\mathbb{Q}\in P$ is specified by $\theta$. There is no misspecification if $\theta \equiv 0$.

We have a 4-D stochastic system to be optimized. The variables $g$ and $\theta$ are the control variables belonging to

$$\mathcal{G} = \{(g_t)_{t\geq 0} | g_t \text{ is -predictable and valued in } \omega \text{ for each } t \geq 0\} \tag{5}$$

and

$$\mathcal{T} = \{(\theta(t,\cdot))_{t\geq 0} | \theta(t,z) \text{ is -predictable and valued in } (-1,+\infty) \text{ for each } t\geq 0\}. \tag{6}$$

Assume $\int_0^\infty \theta^2(t,z)\nu(dz) < +\infty$ a.s. to well-pose (4) and also path-wise unique existence of $X,Y,U,M$ for each $(g,\theta)\in \mathcal{G}\times\mathcal{T}$.

## 2.2 Performance functional

The performance functional $J:[0,T]\times\mathbb{R}_+^4\times\mathcal{G}\times\mathcal{T}\to\mathbb{R}$ is a dynamic risk measure considering the minimization of downstream flood damage (first term) and diversion effort (second term) subject to entropic penalization of the misspecification (third term):

$$\begin{aligned}&J(t,x,y,u,m;g,\theta)\\&=\mathbb{E}_\mathbb{Q}\left[wX_T + w'U_T - \ln M_T | (X_t,Y_t,U_t,M_t)=(x,y,u,m)\right]\\&=\mathbb{E}_\mathbb{P}\left[M_T(wX_T + w'U_T - \ln M_T) | (X_t,Y_t,U_t,M_t)=(x,y,u,m)\right]\end{aligned} \tag{7}$$

where $\mathbb{E}_\mathbb{P}$ (resp., $\mathbb{E}_\mathbb{Q}$) is the expectation under $\mathbb{P}$ (resp., $\mathbb{Q}$) and $w, w' > 0$ are weights as inverses of risk-aversion parameters; larger $w, w'$ represent more risk aversion. The value function is the saddle point

$$\Phi(t,x,y,u,m) = \inf_{g \in \mathcal{G}} \sup_{\theta \in \mathcal{T}} J(t,x,y,u,m;g,\theta). \tag{8}$$

Our goal is to find the optimizer $(g,\theta) = (g^*,\theta^*)$ realizing (8). A practical interest would be $\Phi(0,0,y,0,1)$: the value starting from $(X_0,Y_0,U_0,M_0) = (0,y,0,1)$.

### 2.3 HJBI equation

By a dynamic programming principle, we formally get the HJBI equation:

$$\inf_{g \in \omega} \sup_{\theta:\mathbb{R}_+ \to (-1,+\infty)} L^{g,\theta}\Phi = 0, \quad t \in [0,T), \quad x,y,u,m > 0 \tag{9}$$

with $\Phi(T,x,y,u,m) = m(wx + w'u - \ln m)$ for $x,y,u,m > 0$,

$$L^{g,\theta}\Phi = \frac{\partial \Phi}{\partial t} + (1-g)(y+q)\frac{\partial \Phi}{\partial x} + \frac{g^n}{n}(y+q)\frac{\partial \Phi}{\partial u} - \lambda y \frac{\partial \Phi}{\partial y}$$
$$+ \int_0^\infty \left( \Phi(t,x,y+z,u,m(1+\theta(z))) - \Phi(t,x,y,u,m) - \theta(z)m\frac{\partial \Phi}{\partial m} \right) v(\mathrm{d}z) \tag{10}$$

We can exchange the order of inf and sup in (9) because they are decoupled. With an abuse of notations, the optimizers are expressed as

$$g_t^* = \arg\min_{g \in \omega} \left\{ (1-g)(y+q)\frac{\partial \Phi}{\partial x} + \frac{g^n}{n}(y+q)\frac{\partial \Phi}{\partial u} \right\}, \tag{11}$$

$$\theta^*(t,z) = \arg\max_{\theta:\mathbb{R}_+ \to (-1,+\infty)} \int_0^\infty \left( \Phi(t,x,y+z,u,m(1+\theta(z))) - \theta(z)m\frac{\partial \Phi}{\partial m} \right) v(\mathrm{d}z), \quad z \geq 0, \tag{12}$$

where the right-hand sides are evaluated at $(t,x,y,u,m) = (t,X_t,Y_t,U_t,M_t)$.

## 3. Mathematical analysis

### 3.1 Explicit solution

**Proposition 1** is the main result of this paper.

***Proposition 1***

*(9)-(10) admits a smooth solution*

$$\Phi(t,x,y,u,m) = wxm + w'um - m\ln m + A_t my + B_t m \in C^1([0,T]\times\mathbb{R}_+^4) \qquad (13)$$

*with*

$$A_t = \frac{\hat{w}}{\lambda}\left(1-e^{\lambda(t-T)}\right) \quad \text{and} \quad B_t = \int_t^T\left(q\hat{w} + \int_0^\infty\left(e^{A_t z}-1\right)v(\mathrm{d}z)\right)\mathrm{d}t, \qquad (14)$$

$$\hat{w} = (1-\bar{g})w + n^{-1}\bar{g}^n w' \quad \text{and} \quad \bar{g} = \arg\min_{g\in\omega}\left\{(1-g)w + n^{-1}g^n w'\right\}, \quad \text{if} \quad \left(e^{A_t z}-1\right)v(\mathrm{d}z) \quad \text{is}$$

*integrable in* $(0,\infty)$.

**Proof**

We use a guessed solution technique. Substituting (13) into (10) yields

$$\inf_{g\in\omega}\left\{\frac{\partial\Phi}{\partial t} + (1-g)(y+q)\frac{\partial\Phi}{\partial x} + \frac{g^n}{n}(y+q)\frac{\partial\Phi}{\partial u} - \lambda y\frac{\partial\Phi}{\partial y}\right\}$$

$$= \left(\frac{\mathrm{d}A_t}{\mathrm{d}t}ym + \frac{\mathrm{d}B_t}{\mathrm{d}t}m + (1-\bar{g})(y+q)wm + n^{-1}\bar{g}^n(y+q)w'm - \lambda y A_t m\right) \qquad (15)$$

$$= m\left(\frac{\mathrm{d}A_t}{\mathrm{d}t}y + \frac{\mathrm{d}B_t}{\mathrm{d}t} + \hat{w}(y+q) - \lambda y A_t\right)$$

and

$$\int_0^\infty\left(\Phi(t,x,y+z,u,m(1+\theta(z))) - \Phi(t,x,y,u,m) - m\theta(z)\frac{\partial\Phi}{\partial m}\right)v(\mathrm{d}z)$$
$$= -m\int_0^\infty\left((1+\theta(z))\ln(1+\theta(z)) - z\theta(z)A_t - \theta(z)\right)v(\mathrm{d}z) + m\int_0^\infty zv(\mathrm{d}z)A_t \qquad (16)$$

The optimal $\theta = \theta^*(\cdot)$ should maximize

$$I(\theta) = -\int_0^\infty\left((1+\theta(z))\ln(1+\theta(z)) - A_t z\theta(z) - \theta(z)\right)v(\mathrm{d}z). \qquad (17)$$

Taking the variation of (17) yields

$$\delta I(\theta) = -m \int_0^\infty \left( \ln(1+\theta(z)) - A_t z \right) \delta\theta v(\mathrm{d}z). \tag{18}$$

Therefore,

$$\ln(1+\theta^*(z)) - A_t z = 0 \quad \text{or equivalently} \quad \theta^*(z) = e^{A_t z} - 1(> -1). \tag{19}$$

Then, (9) becomes

$$0 = m \left( \begin{array}{l} \dfrac{\mathrm{d}A_t}{\mathrm{d}t} y + \dfrac{\mathrm{d}B_t}{\mathrm{d}t} + (y+q)\hat{w} - \lambda y A_t \\ - \int_0^\infty \left( (1+\theta^*(z))\ln(1+\theta^*(z)) - z\theta^*(z)A_t - \theta^*(z) \right) v(\mathrm{d}z) + \int_0^\infty z v(\mathrm{d}z) A_t \end{array} \right). \tag{20}$$

Considering (20), we get

$$\dfrac{\mathrm{d}A_t}{\mathrm{d}t} = \lambda A_t - \hat{w}, \quad \dfrac{\mathrm{d}B_t}{\mathrm{d}t} = -q\hat{w} - \int_0^\infty \left( e^{A_t z} - 1 \right) v(\mathrm{d}z), \quad A_T, B_T = 0 \tag{21}$$

and thus (14).

□

The risk-neutral case follows under $\min\{w, w'\} \to +\infty$. **Corollary 1** is a consequence of **Proposition 1** that can be proven in the conventional way [8, Theorem 6.1] owing to the smoothness of $\Phi$.

### *Corollary 1*

*If $\omega = [0,1]$ and $\int_0^\infty (e^{A_t z} - 1) v(\mathrm{d}z) < +\infty$ for $0 \leq t \leq T$, then the optimizers are continuously weight-dependent:*

$$g^* = \begin{cases} 1 & (w \geq w') \\ \left( \dfrac{w}{w'} \right)^{\frac{1}{n-1}} & (w < w') \end{cases} \quad \text{and} \quad \theta^* = e^{A_t z} - 1. \tag{22}$$

The optimal diversion ratio $g^*$ is increasing (resp., decreasing) with respect to the weight $w$ (resp., $w'$), and is simple to implement. By Laeven and Stadjie [17], the worst-case Levy measure is $e^{A_t z} v(\mathrm{d}z)$, meaning that the reference model should be distorted to overestimate both the jump (flood) frequency and size. The coefficients $A_t, B_t$, both are concave, are increasing and decreasing in $t$, respectively. In the long run-limit $T \gg 1$, we get $A_0 \approx \frac{\hat{w}}{\lambda}$, $\frac{B_0}{T} \approx q\hat{w} + \int_0^\infty (e^{A_0 z} - 1) v(\mathrm{d}z)$, and $\Phi(0,0,y,0,1) = A_0 y + B_0$.

*Remark 1*

A maximization problem under uncertainty like a draught mitigation problem can also be analyzed similarly.

### 3.2 Brief application

A critical issue is the well-definiteness of the integral in (14). Recently, we have identified the tempered stable model for the streamflow in a river [5]

$$v(\mathrm{d}z) = \frac{\lambda a}{z^{1+\alpha}} e^{-bz} \mathrm{d}z \quad \text{with} \quad a, b > 0, \ 0 < \alpha < 1. \tag{23}$$

By the asymptotic estimates

$$(e^{A_t z} - 1) v(\mathrm{d}z) \approx O(z^{-\alpha}) \mathrm{d}z, \ 0 < z \ll 1 \tag{24}$$

and

$$(e^{A_t z} - 1) v(\mathrm{d}z) \approx O\left(e^{-(b-A_t)z} z^{-1-\alpha}\right) \mathrm{d}z, \ z \gg 1, \tag{25}$$

the condition $A_t \leq b$ guarantees integrability of $(e^{A_t z} - 1) v(\mathrm{d}z)$. **Proposition 1** then means that $\hat{w}$ (or $\min\{w, w'\}$) or $T$ should be sufficiently small; the decision-maker should be sufficiently close to risk-neutral or the time horizon be sufficiently short. Notice that we may actually require the stronger one $A_t \leq 0.5b$ by the requirement $\int_0^\infty \theta^2(t,z) v(\mathrm{d}z) < +\infty$. Interestingly, these are hard thresholds; for example, we get the locally bounded $\Phi$ if $A_t \leq b$, while we get $\Phi \equiv +\infty$ if $A_t > b$. The problem structure thus suddenly changes across the threshold. Therefore, the presented model is robust for decision-makers close to risk-neutral.

We conclude this section with a backward SDE (BSDE)-based alternative of the HJBI equation. It is well-known that HJBI equations and BSDEs are closely linked with each other [18]. In our case, setting $V_t = \Phi(t, X_t, Y_t, U_t, M_t)$ with $g = g^*$ yields

$$\mathrm{d}V_t = -\int_0^{+\infty} \left(e^{Z(t,z)} - Z(t,z) - 1\right) v(\mathrm{d}z)\mathrm{d}t + \int_0^{+\infty} Z(t,z) \tilde{N}(\mathrm{d}t, \mathrm{d}z), \quad 0 \leq t < T \qquad (26)$$

with a $\mathbb{F}$-predictable $Z(t,z) = A_t z$. This kind of alternative representation provides a dual view of both the presented and further extended problems in future.

## 4. Conclusion

We presented an exactly-solvable model of river water management using a dynamic risk measure based on the entropic penalization. The explicit nature of the solution to the HJBI equation clearly demonstrated its parameter dependence and the threshold nature of the risk-aversion behavior.

We used an entropic penalty, but non-entropic ones can also be used with slight modifications [19]. The BSDE representation opens a door to numerically handle cases with a larger number of variables. Accurate approximations of $Z$ and $\tilde{N}$ under infinite-activities Lévy processes would be an important task.

**Acknowledgements**

Kurita Water and Environment Foundation 19B018 and 20K004 support this research.

**References**

[1] Lundström, N. L., Olofsson, M., & Önskog, T. (2020). Management strategies for run-of-river hydropower plants-an optimal switching approach. arXiv preprint arXiv:2009.10554.

[2] Botter, G., Zanardo, S., Porporato, A., Rodriguez-Iturbe, I., & Rinaldo, A. (2008). Ecohydrological model of flow duration curves and annual minima. Water Resources Research, 44(8), W08418. https://doi.org/10.1029/2008WR006814

[3] Ramirez, J. M., & Constantinescu, C. (2020). Dynamics of drainage under stochastic rainfall in river networks. Stochastics and Dynamics, 20(03), 2050042. https://doi.org/10.1142/S0219493720500422

[4] Bianchi, M. L., Rachev, S. T., & Fabozzi, F. J. (2017). Tempered stable Ornstein–Uhlenbeck processes: A practical view. Communications in Statistics-Simulation and Computation, 46(1), 423-445. https://doi.org/10.1080/03610918.2014.966834

[5] Yoshioka, H., & Yoshioka, Y. (2020). Tempered stable Ornstein–Uhlenbeck model for river discharge time series with its application to dissolved silicon load analysis, International Conference on Water Security and Management, Dec. 15-18, Tokyo, Japan, Online. Accepted on 2020/10/15. Full paper 10pp.

[6] Yoshioka, H., Yoshioka, Y. (2020). Regime switching constrained viscosity solutions approach for controlling dam-reservoir systems. Computers and Mathematics with

Applications, 80(9), 2057-2072. https://doi.org/10.1016/j.camwa.2020.09.005

[7] Zakaria, A., Ismail, F. B., Lipu, M. H., & Hannan, M. A. (2020). Uncertainty models for stochastic optimization in renewable energy applications. Renewable Energy, 145, 1543-1571. https://doi.org/10.1016/j.renene.2019.07.081

[8] Øksendal, B., & Sulem, A. (2019). *Applied Stochastic Control of Jump Diffusions*. Springer, Cham.

[9] Föllmer, H., & Schied, A. (2002). Convex measures of risk and trading constraints. Finance and stochastics, 6(4), 429-447. https://doi.org/10.1007/s007800200072

[10] Herman, J. D., Quinn, J. D., Steinschneider, S., Giuliani, M., & Fletcher, S. (2020). Climate adaptation as a control problem: Review and perspectives on dynamic water resources planning under uncertainty. Water Resources Research, 56(2), e24389. https://doi.org/10.1029/2019WR025502

[11] Wellen, C., Van Cappellen, P., Gospodyn, L., Thomas, J. L., & Mohamed, M. N. (2020). An analysis of the sample size requirements for acceptable statistical power in water quality monitoring for improvement detection. Ecological Indicators, 118, 106684. https://doi.org/10.1016/j.ecolind.2020.106684

[12] Wei, C., & Luo, C. (2020). A differential game design of watershed pollution management under ecological compensation criterion. Journal of Cleaner Production, 274, 122320. https://doi.org/10.1016/j.jclepro.2020.122320

[13] Jiang, K., You, D., Li, Z., & Shi, S. (2019). A differential game approach to dynamic optimal control strategies for watershed pollution across regional boundaries under eco-compensation criterion. Ecological Indicators, 105, 229-241. https://doi.org/10.1016/j.ecolind.2019.05.065

[14] Komaee, A. (2020). An inverse optimal approach to design of feedback control:

Exploring analytical solutions for the Hamilton-Jacobi-Bellman equation. Optimal Control Applications and Methods. https://doi.org/10.1002/oca.2686

[15] Yoshioka, H., & Yoshioka, Y. (2019). Modeling stochastic operation of reservoir under ambiguity with an emphasis on river management. Optimal Control Applications and Methods, 40(4), 764-790. https://doi.org/10.1002/oca.2510

[16] Cont R., & Tankov, P. (2003). *Financial Modelling with Jump Processes*. CRC press, Boca Raton, London, New York, Washington, D.C..

[17] Laeven, R. J., & Stadje, M. (2014). Robust portfolio choice and indifference valuation. Mathematics of Operations Research, 39(4), 1109-1141. https://doi.org/10.1287/moor.2014.0646

[18] Delong, Ł. (2013). *Backward Stochastic Differential Equations with Jumps and Their Actuarial and Financial Applications*. Springer, London.

[19] Faidi, W., Matoussi, A., & Mnif, M. (2017). Optimal stochastic control problem under model uncertainty with nonentropy penalty. International Journal of Theoretical and Applied Finance, 20(03), 1750015. https://doi.org/10.1142/S0219024917500157